\documentclass[preprint,12pt]{elsarticle}

\usepackage{graphicx}
\usepackage{subfigure}
\usepackage{amssymb}
\usepackage[english]{babel}
\usepackage[latin1]{inputenc}
\usepackage{amsmath}
\usepackage{color}

\usepackage{hyphenat}

\newtheorem{theorem}{Theorem}

\newtheorem{remark}{Remark}

\begin{document}

\begin{frontmatter}


% use the tnoteref command within \title for footnotes;
% use the tnotetext command for theassociated footnote;
% use the fnref command within \author or \address for footnotes;
% use the fntext command for theassociated footnote;
% use the corref command within \author for corresponding author footnotes;
% use the cortext command for theassociated footnote;
% \author{Name\corref{cor1}\fnref{label2}}
% \address{Address\fnref{label3}}

\title{Optimal management of open-channel raceway ponds for cultivation of algal biomass intended for bioenergy production}


\author[vigo]{L.J. Alvarez-V\'azquez}
\ead{lino@dma.uvigo.es}
\address[vigo]{Departamento de Matem\'atica Aplicada II, Universidade de Vigo. \\ E.E.~Telecomunicaci\'on, 36310 Vigo. Spain}
\author[vigo]{A. Mart\'inez}
\ead{aurea@dma.uvigo.es}
%\author[santiago]{C. Rodr\'iguez}
%\ead{carmen.rodriguez@usc.es}
%\address[santiago]{Departamento de Matem\'atica Aplicada, Universidade de
%Santiago de Compostela. \\ Fac. Matem\'aticas, 15782 Santiago.  Spain}
\author[lugo]{M.E.~V\'azquez-M\'endez}
\ead{miguelernesto.vazquez@usc.es}
\address[lugo]{Departamento de Matem\'atica Aplicada, Universidade de Santiago de Compostela. \\ E.P.S.  Enxe\~ner\'{\i}a, 27002 Lugo.  Spain}

\begin{abstract}
In this work we present a novel methodology to deal with the optimal perfor\-mance of raceways (open-channel ponds where the circulating wastewater, during its purification process, is used to grow algae that will be used as a source for the production of bioenergy).  The maximization of algal producti\-vity is addressed here within an optimal control framework for partial differen\-tial equations.  Thus, after introducing a rigorously detailed mathematical formulation of the real-world control problem,  we prove the existence of optimal solutions,  we propose a numerical algorithm for its computational resolution and,  finally,  we show some results for the numerical optimization of a realistic case.
\end{abstract}

\begin{keyword}
Optimal control \sep Mathematical modelling \sep Simulation-based optimization \sep Open raceway pond \sep Algaculture \sep Bioenergy production 
\MSC 35K57 \sep 35K55 \sep 49J20 \sep 74S05 \sep 65K05
% MSC codes here, in the form: \MSC code \sep code
% or \MSC[2008] code \sep code (2000 is the default)
\end{keyword}

\end{frontmatter}

\section{Introduction}

Treatment of wastewater -usually from domestic,  industrial,  agricultural or livestock origin- by means of microalgae-based technologies is nowadays an effective solution that allows, in addition, the recovery of resources and materials which can be used, for instance, for the production of bioenergy (mainly biodiesel, but also bioalcohol,  methane,  biohydrogen, etc.),  since algal biomass from wastewater treatment is rich in lipides \cite{Sha} that represent a promising alternative source of oil for the production of bioenergy.

There exists a wide range of algae cultivation systems: open, closed and hybrid ones (interested readers can find an exhaustive overview of the topic in recent review paper \cite{Bha}).  However,  among the most commonly used systems of algae cultivation,  open cultivation systems stand out,  in particular the artificially designed bodies of water known as raceway ponds (as opposed to formerly used natural ponds such as lakes or lagoons).  Raceways constitute the cheapest and more advantageous option among all the simple choices due to their easy construction process, low maintenance costs, simple operation,  and low energy consumption.  Raceway ponds are a very effective way of algae growth and harvesting: other alternatives, such as photobioreactors or algal turf scrubbers,  can be more productive and easier to control,  and less dependent on climate (light intensity,  temperature,  atmospheric CO$_2$ and O$_2$ levels...),  but they are dramatically more expensive.  

A raceway cultivation system is a shallow artificially engineered pond formed by an open channel in the shape of an oval (the appellation {\it raceway} comes from its similarity to an automotive raceway circuit) equipped with a rotating paddlewheel intended to avoid algal sedimentation and promote water recirculation for a more productive algal culture.  Raceway ponds have been used since 1950s.  In order to employ wastewater as a resource for energy production,  algae play a fundamental role, due both to their ability to remove nutrients from untreated water and to their capacity to accumulate lipids for a simultaneous production of bioenergy.  Raceways allow an effective algae cultivation using different types of wastewaters, where -once the growth of the algal mass has reached the desired level- this algal biomass must be recovered by various types of mechanical,  chemical,  or biological harvesting methods,  such as filtration, sedimentation, flocculation or centrifugation.  Harvested algae biomass can be then used as feedstock for bioenergy production (and also for alternative biochemical production).

\begin{figure}[t]
\centering
\includegraphics[width=11cm]{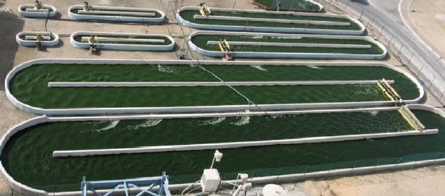}
\caption{Examples of real-world raceway ponds of different sizes.} \label{fig1}
\end{figure}

Once the species of microalgae to be cultivated has been chosen \cite{Aze} (basi\-cally depending on the characteristics of influent wastewater),  one of the main external factors -easily managed by raceway stakeholders- that influences algal productivity is related to water velocity inside the pond.  This velocity,  which can be directly controlled by the position and rotational speed of the turning paddle wheel, greatly conditions algal growth, mainly reducing dead zones,  preventing sedimentation of algae,  and promoting mixing and aeration of water. 

The study, both analytical and numerical, of the modelling of algae growth in raceway ponds for bioenergetic purposes is a subject that has been widely addressed during the last decades.  For the case of the mathematical modelling of the process, readers can consult the recent, extensive critical review \cite{Sho}, where more than 320 published models (from both articles and conference proceedings) are presented and compared.  For the use of computa\-tional fluid dynamics in order to simulate the performance of ponds a detailed review can be found, for instance,  in the recent paper \cite{Pir}, where the most frequently used models (CFX, Fluent, COMSOL\dots) are discussed.  
The optimization and design of raceway ponds has also been extensively analyzed, but mainly from the point of view of economic productivity,  pond geometry or energy saving.  Very different approaches can be found, for example, in the interesting works \cite{Ray, Rog, Hre, Hua, Yad, Ban} and the references therein.
However, the optimal management of the rotating paddlewheel or the volume of water have not been, as far as we know, as thoroughly studied. We can mention here the works of Pandey-Premalatha \cite{Pan}, Chen {\it et al.} \cite{Che} or Ali {\it et al.} \cite{Ali}, where the speed -and in some cases the position- of the paddlewheel is optimized, but essentially by comparing case studies within a statistical framework. The use of the techniques of optimal control of partial differential equations for the simultaneous optimization of water volume and speed of the rotating paddlewheel has so far remained completely unaddressed within the mathematical litera\-ture.  A much simplified alternative approach to some related control problems,  based on only time-dependent ordinary differential equations, can be seen in \cite{Ber2} and therein references.

In this work,  we introduce a rigorously detailed mathematical formulation of the optimal control problem,  we demonstrate its solvability -proposing a full numerical algorithm for its resolution- and, finally,  we show a few  computational results related to the numerical optimization of the problem.

\section{Mathematical setting of the problem}

\subsection{The state systems}

\begin{figure}[t]
\centering
\includegraphics[width=11cm]{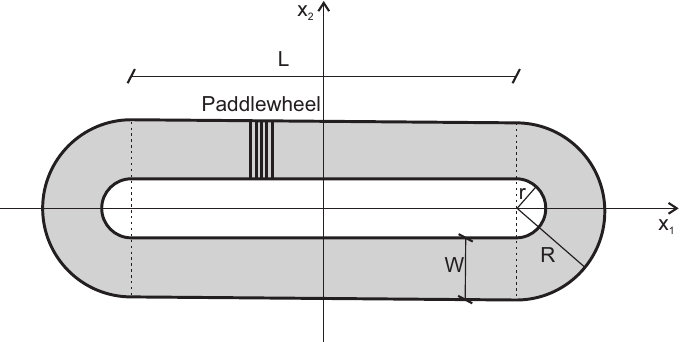}
\caption{Schematic drawing of the raceway ground plan $G$, showing the two straight channels of length $L$ and width $W$,  and the two semicircular channels of radii $r$ and $R=r+W$.  A possible location for the paddlewheel is also shown.} \label{figG}
\end{figure}
 
We consider a moving liquid domain $\Omega(t) \subset \mathbb{R}^3$, for each $t$ in a time interval $I = (0,T)$,  representing the open raceway pond occupied by shallow waters:
$$\Omega(t) = \{ (x_1,x_2,x_3) \in \mathbb{R}^3 \ / \  (x_1,x_2) \in G,  \ 0 < x_3 < \eta(x_1,x_2,t) \} $$
where $G\subset \mathbb{R}^2$ is a fixed domain, with a regular boundary $\partial G$,  defining the fishway ground plan (as given in Fig.~\ref{figG}) and $\eta$ represents the height of the water column.
Denoting by $\Gamma(t)$ the boundary of $\Omega(t)$ (assumed to be smooth enough), we suppose this boundary $\Gamma(t)$ divided into two parts $\Gamma(t) = \Gamma_1(t) \cup \Gamma_2(t)$, where $\Gamma_1(t)$ corresponds to the bottom (i.e., $(x_1,x_2) \in G,  \, x_3=0$) and the lateral walls (i.e., $(x_1,x_2) \in\partial G,  \, 0 < x_3 < \eta(x_1,x_2,t)$),  and $\Gamma_2(t)$ corresponds to the top free surface (i.e., $(x_1,x_2) \in G,  \, x_3=\eta(x_1,x_2,t)$).  Finally,  we denote by $\vec{n}(x_1,x_2,x_3,t)$ the unit outward normal vector to the boundary $\Gamma(t)$, and by $\vec{t}(x_1,x_2,x_3,t)$ a tangential unit vector to the boundary.

\subsubsection{The hydrodynamic model} \label{sub1}

Our first set of state variables includes velocity $\vec{v}(x,t)=(v_1,v_2,v_3)$ of the liquid at time $t \in I$ and at point $x \equiv (x_1,x_2,x_3) \in \Omega(t)$ and pressure $p(x,t)$ (decomposed as a hydrostatic component plus a hydrodynamic one),  given by the classical Navier-Stokes equations for incompressible flows with free surface (where the effect of surface tension is neglected).  So, for $Q=\cup_{t\in I} \, \Omega(t) \times \{t\}$,  $\Sigma_1=\cup_{t\in I} \, \Gamma_1(t) \times \{t\}$ and $\Sigma_2=\cup_{t\in I} \, \Gamma_2(t) \times \{t\}$,  we have:
\begin{equation} \label{v}
\left\{ \begin{array}{l}
\displaystyle{\frac{\partial \vec{v} }{\partial t}+(\vec{v} \cdot \nabla) \vec{v} -\mu \Delta \vec{v}  + %\frac{1}{\rho}
\nabla p} = %\vec{g} + 
\vec{F} \ \ \ {\rm in } \ Q ,
\vspace*{0.15cm} \\
\displaystyle{\nabla \cdot \vec{v} = 0} \ \ \ {\rm in } \ Q ,
\vspace*{0.15cm} \\
\displaystyle{\vec{v} \cdot \vec{n}}= 0 \ \ \ {\rm on } \ \Sigma_1 ,
\vspace*{0.15cm} \\
\displaystyle{(\mu \nabla\vec{v} - p {\mathcal I}) \vec{n} \cdot \vec{t} = 0} \ \ \ {\rm on } \ \Sigma_1 ,
\vspace*{0.15cm} \\
\displaystyle{(\mu \nabla\vec{v} - p {\mathcal I}) \vec{n} = \vec{0}} \ \ \ {\rm on } \ \Sigma_2 ,
\vspace*{0.15cm} \\
\displaystyle{\vec{v}(x,0)} = \vec{v}_0(x) \ \ \ {\rm in } \ \Omega(0) ,
\end{array}\right.
\end{equation}
where $\mu$ is the dynamic viscosity coefficient,  ${\mathcal I}$ is the identity matrix,  %$\rho$ is the density,  $\vec{g}=(0,0,-g)$ stands for gravity force,  
and the second-member forcing term $\vec{F}(x,t)=(F_1, F_2,F_3)$ represents the effect of the rotating paddlewheel with its axis centered at point $(x_1^0,x_2^0,x_3^0)$, with paddles of length $\rho$,  and with angular speed $\omega$ (see full details, for instance, in \cite{Ber}), that is,  for a force magnitude $F$:
\begin{equation} \label{F}
\begin{array}{l}
F_1(x,t) = F \omega^2 \cos(\omega t) [(x_1-x_1^0)^2+(x_3-x_3^0)^2],  \\
F_2(x,t) = 0 ,  \\
F_3(x,t) = F \omega^2 \sin(\omega t) [(x_1-x_1^0)^2+(x_3-x_3^0)^2] ,  
\end{array}
\end{equation}
in the region of influence of the paddles $R(t)$ (further details about its definition will be given in below Section \ref{S4}, but $R(t)$ represents the intersection of the cylinder $R = \{ (x_1, x_2,x_3) \in \mathbb{R}^3 \ /  \  x_2^0 - W/2 \leq x_2 \leq x_2^0 + W/2 , \  (x_1-x_1^0)^2 + (x_3-x_3^0)^2 \leq \rho^2 \} $ with $\Omega(t)$),  and $ \vec{F} $ is null in the rest of the raceway.

Navier-Stokes equations must be completed with the following kinematic condition at the free surface:
\begin{equation} \label{K}
%\displaystyle{\frac{\partial\eta}{\partial t} + v_1 \frac{\partial\eta}{\partial x_1} + v_2 \frac{\partial\eta}{\partial x_2} = v_3} \ \ \ {\rm on } \ \Gamma_2(t), \ t \in I ,
\displaystyle{\frac{\partial\eta}{\partial t} + v_{1|x_3=\eta} \frac{\partial\eta}{\partial x_1} + v_{2|x_3=\eta} \frac{\partial\eta}{\partial x_2} = v_{3|x_3=\eta} } \ \ \ {\rm in } \ G \times I ,
\end{equation}
(representing the fact that fluid particles at free surface remain at the surface, assuring the conservation of the liquid quantity),
and also the initial condition:
\begin{equation} \label{K2}
\displaystyle{\eta(x_1,x_2,0) = \eta_0(x_1,x_2) }  \ \ \ {\rm in } \ G.
\end{equation}

In our real-world case,  we will assume that, at initial time $t=0$,  the pond presents a fixed constant height of water $H>0$ (that is,  $\eta_0=H$ or, equivalently,  $\Omega(0)= G \times (0,H)$),  where water is initially at rest (that is,  $\vec{v}_0=\vec{0}$).  Moreover, due to the characteristics of our problem, it is expected that, after a short transition period, the hydrodynamic system will reach a stable periodic flow regime \cite{Nik}.

\subsubsection{The biological model}

We also consider the state variables corresponding to algal concentration $A(x,t)$,  $PO_4$ concentration $P_1(x,t)$,  non-assimilable $P$ concentration $P_2(x,t)$,  $NO_3$ concentration $N_1(x,t)$,  non-assimilable $N$ concentration $N_2(x,t)$,  $NH_4$ concentration $N_3(x,t)$,  organic load $D(x,t)$ and dissolved oxygen $O(x,t)$,  given by the following coupled nonlinear system of convection-diffusion-reac\-tion equations with Monod kinetics \cite{Al1, Al2},  posed in $Q$,  and with liquid velocity $\vec{v}$ obtained from previous state system (\ref{v}):

\begin{equation} \label{A} %\nonumber
\hspace*{-1.35cm} \left\{ \begin{array}{l}
\displaystyle{\frac{\partial A }{\partial t}+\vec{v} \cdot \nabla A -\mu_A\Delta A  = (L \, \frac{P_1}{K_P+P_1} \,\frac{N_1+N_3}{K_N+N_1+N_3} - (\gamma + \beta))  A } ,
\vspace*{0.15cm} \\
\displaystyle{\frac{\partial P_1 }{\partial t}+\vec{v} \cdot \nabla P_1 -\mu_P\Delta P_1  = C_P (\delta_1 (\gamma + \beta) - L \,\frac{P_1}{K_P+P_1}\,\frac{N_1+N_3}{K_N+N_1+N_3})A }
\vspace*{0.15cm} \\ 
\qquad \qquad + \kappa_1 P_2 ,
\vspace*{0.15cm} \\
\displaystyle{\frac{\partial P_2 }{\partial t}+\vec{v} \cdot \nabla P_2 -\mu_P\Delta P_2  = C_P (1-\delta_1) (\gamma + \beta) A } - \kappa_1 P_2 - W P_2  ,
\vspace*{0.15cm} \\
\displaystyle{\frac{\partial N_1 }{\partial t}+\vec{v} \cdot \nabla N_1 -\mu_N\Delta N_1  = - C_N L \, \frac{P_1}{K_P+P_1}\,\frac{N_1}{K_N+N_1+N_3} A + \kappa_2 N_3}  ,
\vspace*{0.15cm} \\
\displaystyle{\frac{\partial N_2}{\partial t}+\vec{v} \cdot \nabla N_2 -\mu_N\Delta N_2  = C_N (1-\delta_2) (\gamma + \beta) A - \kappa_3 N_2 - W N_2 } ,
\vspace*{0.15cm} \\
\displaystyle{\frac{\partial N_3 }{\partial t}+\vec{v} \cdot \nabla N_3 -\mu_N\Delta N_3  = C_N (\delta_2 (\gamma + \beta) - L \, \frac{P_1}{K_P+P_1}\,\frac{N_3}{K_N+N_1+N_3}) A }
\vspace*{0.15cm} \\ 
\qquad \qquad + \kappa_3 N_2 - \kappa_2 N_3  ,
\vspace*{0.15cm} \\
\displaystyle{\frac{\partial D }{\partial t}+\vec{v} \cdot \nabla D -\mu_D\Delta D  = \phi \gamma A - \kappa_4 \Theta^{\theta(t)-\theta_0} D - W D } ,
\vspace*{0.15cm} \\
\displaystyle{\frac{\partial O }{\partial t}+\vec{v} \cdot \nabla O -\mu_O\Delta O  = \phi ( L \,\frac{P_1}{K_P+P_1} \, \frac{N_1+N_3}{K_N+N_1+N_3} - \beta ) A } 
\vspace*{0.15cm} \\ 
\qquad \qquad - \nu \kappa_2 N_3 - \kappa_4 \Theta^{\theta(t)-\theta_0} D + \kappa_3 \Theta^{\theta(t)-\theta_0} (C_s - O) - B ,
\end{array}\right.
\end{equation}
where the light rays effect $L$ on algae is given by the expression:
$$L(x,t) = \mu_{max} \Theta^{\theta(t)-\theta_0} i(t) \, e^{-(\Phi_1 + \Phi_2 A) x_3}$$
with $\mu_{max}$ the maximum specific growth rate,  $\Theta$ the thermic regeneration coefficient,  $\theta(t)$ the temperature,  $\theta_0$ a reference temperature, $i(t)$ the incident light intensity,  and $\Phi_1$,  $\Phi_2$ the coefficients for light attenuation due to depth and algal mass.  The other parameters in state system (\ref{A}) are associated with the various physical-chemical-biological phenomena that combine the different species.  These constant parameters are, specifically,  $\mu_A$,  $\mu_N$, $\mu_P$,  $\mu_D$ and $\mu_O$ (corresponding to the diffusion coefficients),  $\gamma$ (the algal death rate),  $\beta$ (the algal respiration rate),  $K_N$ and $K_P$ (the half-saturation constant for nitrogen and phosphorus,  respectively),  $C_N$ and $C_P$ (representing the stoichiometric relations for nitrogen and phosphorus,  respectively),  $\delta_1$ (the proportion of assimilable $P$ in dead algae),  $\kappa_1$(the transformation rate of non-assimilable $P$ into $PO_4$),  $W$ (the sedimentation velocity,)  $\kappa_2$ (the kinetic nitrification constant),  $\delta_2$ (the proportion of assimilable $N$ in dead algae),  $\kappa_3$ (the transformation rate of non-assimilable $N$ into $NO_3$),  $\phi$ (the oxygen quantity produced by photosynthesis),  $\kappa_4$ (the kinetic degradation constant of organic load),  $\nu$ (the oxygen consumed in nitrification),  $C_s$ (the saturation concentration of oxygen),  and $B$ (the benthic oxygen demand).

System (\ref{A}) must be completed with null Neumann boundary conditions on $\Sigma = \Sigma_1 \cup \Sigma_2$ for all variables (corresponding to no-flux conditions), and with given initial conditions $A_0$,  $P_{1,0}$,  $P_{2,0}$,  $N_{1,0}$,  $N_{2,0}$,  $N_{3,0}$,  $D_0$ and $O_0$ bounded in $\Omega(0)$.

\subsection{The optimal control problem}

In this study we are interested in finding the optimal initial height $H$ of water and the optimal rotational speed $\omega$ of the paddlewheel such that the production of algal biomass at the final time of the process is maximal,  that is,  for instance, maximizing a function of the form
$ \int_{\Omega(T)} A(x,T) \, dx $.

Moreover,  for a proper operation of the raceway,  we need to secure a high enough global velocity of water (keeping algae in suspension,  assuring an effective distribution of nutrients,  and avoiding the presence of dark, deep regions at the bottom of the pond),  and we also have to guarantee that the dissolved oxygen concentration remains always over a critical threshold (essential for algal growth).  So, we impose the following state constraints:
\begin{eqnarray} 
&& \int_0^T \int_{\Omega(t)} \| \vec{v}(x,t) \| \, dx \, dt \geq C_1 \geq 0 , \label {c1} \\
&& \min_{t \in I} \int_{\Omega(t)} O(x,t) \, dx \geq C_2 \geq 0 .  \label{c2}
\end{eqnarray}

Thus, taking into account the fact that the maximum of a function corresponds to the minimum of its negative,  %and incorporating above state constraints as penalty terms (with associated weights $M_1,  M_2 >0$),  
we are finally led to solve the optimal control problem $(\mathcal{P})$: finding the optimal values for $H$ and $\omega$ (both subject to appropriate geometric and technological constraints $0 < H_{min} \leq H \leq H_{max}$, $0 < \omega_{min} \leq \omega \leq \omega_{max}$) that satisfy the state constraints (\ref {c1}) and (\ref {c2}),  and minimize the cost functional:
\begin{eqnarray} 
&& J(H,\omega) = - \int_{\Omega(T)} A(x,T) \, dx %+ M_1 \max \left\{ C_1 - \int_0^T \int_{\Omega(t)} \| \vec{v}(x,t) \| \, dx \, dt , 0 \right\} \nonumber \\
%&& \qquad \qquad + M_2 \max \left\{ C_2 - \int_0^T \int_{\Omega(t)} O(x,t) \, dx \, dt , 0 \right\}  
\label{J}
\end{eqnarray}
%or (if the state variable $A$ is not regular enough,  for instance, not continuous at $t=T$) the alternative cost functional:
%\begin{equation} \label{JJ}
%J(H,\omega) = - \int_0^T \int_{\Omega(t)} A(x,t) \, dx \, dt ,
%\end{equation}
where the control variables $H$ and $\omega$ enter the cost function via the initial configuration $\Omega(0)$ and the second member of the system (\ref{v}) corresponding to state variable $\vec{v}$, respectively.  A related optimal control problem,  also involving algal growth in a moving domain,  has been previously analyzed by the authors in \cite{Al3}.

\begin{remark}\label{R1}
Regarding the regularity of the state variables,  we will prove in next section that all of them have enough regularity,  so that all integrals in above optimal control problem $(\mathcal{P})$ make sense.
However,  if any state variable is not smooth enough (for instance,  $A$ is not continuous at $t=T$),  the alternative expression $J = - \int_0^T \int_{\Omega(t)} A(x,t) \, dx \, dt$
could be used in the cost functional instead of $J = - \int_{\Omega(T)} A(x,T) \, dx$,  as given in (\ref{J}).  Similar reasonings can be also employed for suitable alternative formulations of the state constraints (\ref {c1}) and (\ref {c2}).
\end{remark}

\section {Analytical study of the problem}

Our main goal in this Section is to demonstrate the existence of solution of above optimal control problem $(\mathcal{P})$.  For a simpler presentation of our results,  we will define the reference domain $\hat\Omega = G \times (0,1)$,  so that,  for any water height $H \in [H_{min}, H_{max}]$,  initial domain $\Omega(0) = \tau(\hat\Omega)$, where $\tau$ is the smooth mapping given by:
$$\tau : \zeta \equiv (\zeta_1,\zeta_2,\zeta_3) \in \hat\Omega \to \xi = \tau(\zeta_1,\zeta_2,\zeta_3) = (\zeta_1,\zeta_2,H \zeta_3) = \mathcal{H} \zeta \in \Omega(0) , $$
where $\mathcal{H}$ is the diagonal matrix $\mathcal{H}= {\rm diag}(1,1,H)$.  Trivially,  $\nabla \tau = \mathcal{H}$, and $\det(\mathcal{H}) = H > 0$.  Thus,  matrix $\mathcal{H}$ is invertible.

The strategy for the proof of this existence result involves three steps.  The first and the second ones are related to obtaining several {\it a priori} estimates for the solutions of the hydrodynamic system (\ref{v}) and the biological system (\ref{A}), respectively. The third step focuses on the determination of the optimal solution {\it via} minimizing sequences techniques.

\subsection{Estimates for the hydrodynamic model}

We consider the following hydrodynamic problem,  defined for $x \in \Omega(t)$,  $ t \in (0, \infty)$:
\begin{equation} \label{vv}
\left\{ \begin{array}{l}
\displaystyle{\frac{\partial \vec{v} }{\partial t}+(\vec{v} \cdot \nabla) \vec{v} -\mu \Delta \vec{v}  + \nabla p} =  \vec{f} ,
\vspace*{0.15cm} \\
\displaystyle{\nabla \cdot \vec{v} = 0} ,
\end{array}\right.
\end{equation}
with initial condition $\vec{v}(x,0) = \vec{0}$ in $\Omega(0)$,  and boundary conditions as in (\ref{v}).

From a theoretical viewpoint, incompressible Navier-Stokes equations with free surface have been the subject of several papers, although in most of them the authors consider periodic boundary conditions and/or infinite spatial domains,  which unfortu\-nately is not applicable to our case.  So, for instance,  in \cite{exi1, exi2, exi3} some local existence results are demonstrated for the problem in unbounded domains.  For the case of bounded domains,  some existence results have been obtained under the assumption of small data, for example, in \cite{exi4-},  \cite{exi4} or \cite{exi5}.

\begin{figure}[t]
\centering
\includegraphics[width=9cm]{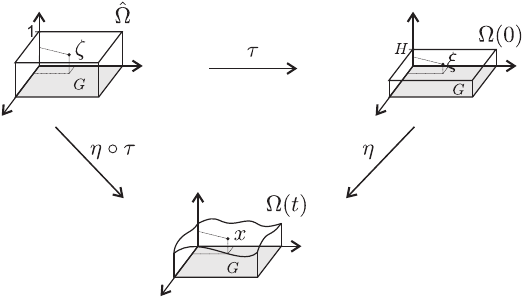}
\caption{Relations between reference coordinates $\zeta\in\hat\Omega$, Lagrange coordinates $\xi\in\Omega(0)$, and Euler coordinates $x\in\Omega(t)$.} \label{fig3}
\end{figure}

In our particular case,  using the seminal results of Solonnikov \cite{exi4} we know that, if $\vec{f}$ is bounded (that is,  if $\int_0^{\infty} \sup_x \| \vec{f}(x,t) \| dt < \infty$),  then the problem (\ref{vv}) has a unique solution $q \in W_p^{1,0}(\hat\Omega\times(0,T_1))$,  $\vec{u} \in [W_p^{2,1}(\hat\Omega\times(0,T_1))]^3$,  for $p>3$,  with $T_1(\sup \| \vec{f} \|) \to \infty$ as $\sup \| \vec{f} \| \to 0$,  in reference coordinates $\zeta \in \hat\Omega$,  that is:
$$\vec{u}(\zeta,t) = \vec{v}(\eta(\tau(\zeta),t),t) = \vec{v}(x,t),  \qquad q(\zeta,t) = p(\eta(\tau(\zeta),t),t) = p(x,t), $$
where Lagrange coordinates $\xi\equiv (\xi_1,\xi_2,\xi_3) = \tau(\zeta) \in \Omega(0)$ are related to Euler coordinates $x=\eta(\xi,t) \in \Omega(t)$ by the following ordinary differential equation:
\begin{equation} \label{eta}
\left\{ \begin{array}{l}
\displaystyle{\frac{\partial \eta(\xi,t) }{\partial t} = \vec{v}(x,t)} ,
\vspace*{0.15cm} \\
\displaystyle{\eta(\xi,0) = \xi} ,
\end{array}\right.
\end{equation}
for $\xi \in  \Omega(0)$, $t \in (0, \infty)$.

\begin{remark}\label{R2}
We recall here that functional space $W_p^{2,1}(\hat\Omega\times(0,T))$ corresponds to all functions $g(\zeta,t)$ such that $g$,  $\frac{\partial g}{\partial t}$,  $ \frac{\partial g}{\partial \zeta_i}$ and $ \frac{\partial^2 g}{\partial \zeta_i \partial \zeta_j}$ belong to $L^{p}(\hat\Omega\times(0,T))$, that is,  $ W_p^{2,1}(\hat\Omega\times(0,T)) = L^{p}((0,T);W^{2,p}(\hat\Omega)) \cap W^{1,p}((0,T);L^{p}(\hat\Omega))$.
In a similar manner,  space $W_p^{1,0}(\hat\Omega\times(0,T))$ denotes all functions $g(\zeta,t)$ such that $g$ and $ \frac{\partial g}{\partial \zeta_i}$ lie in $L^{p}(\hat\Omega\times(0,T))$ or,  equivalently,  $W_p^{1,0}(\hat\Omega\times(0,T)) = L^{p}((0,T);W^{1,p}(\hat\Omega))$.
\end{remark}

\begin{remark}\label{R3}
For our particular case, with $\vec{f} = \vec{F} $ (expanded by zero when $t > T$),  its boundedness is trivial: we have that,  for any $t$ and for any $x \in R(t)$,  $\| \vec{F}(x,t) \| \leq F \omega^2 \rho^2$, and consequently $\sup \| \vec{F} \|  \leq F \omega^2 \rho^2 $. Then, 
$$\int_0^{\infty} \sup_x \| \vec{F}(x,t) \| dt = \int_0^{T} \sup_x \| \vec{F}(x,t) \| dt \leq T F \omega^2 \rho^2 < \infty.$$
On the other hand,  since $\sup \| \vec{F} \|  \leq F \omega^2 \rho^2 $, 
taking $F$ and/or $\rho$ small enough, we obtain the existence result of $q$ and $\vec{u}$ for $T_1=T$, that is,  $q \in W_p^{1,0}(\hat\Omega\times(0,T))$,  $\vec{u} \in [W_p^{2,1}(\hat\Omega\times(0,T))]^3$,  for $p>3$.  Moreover,  since $R(t) \subset R$,  for all $t$,  we also have that the volumetric measure ${\rm meas}(R(t)) \leq {\rm meas}(R) = \pi W \rho^2$.  
\end{remark}

We define the matrix ${\mathcal A}(\xi,t)$ as the cofactor matrix of $\nabla_{\xi} \,\eta$,  where $\nabla_{\xi} \equiv (\frac{\partial}{\partial \xi_1}, \frac{\partial}{\partial \xi_2}, \frac{\partial}{\partial \xi_3})$ represents the gradient with respect to Lagrange coordinates $\xi$.  Moreover, since ${\nabla_x\cdot\vec{v} = 0}$,  where $\nabla_x \equiv (\frac{\partial}{\partial x_1}, \frac{\partial}{\partial x_2}, \frac{\partial}{\partial x_3})$ denotes the gradient with respect to Euler coordinates $x$,  we have that $\det({\mathcal A})=1$,  and, conse\-quently,  the invertible matrix ${\mathcal A}$ can be alternatively defined by ${\mathcal A}^T =(\nabla_{\xi} \,\eta)^{-1}$ \cite{exi4}.  
Then, as can be also seen in \cite{exi4},  the gradient with respect to Lagrange coordinates $\xi$ verifies $\nabla_{\xi} = {\mathcal A} \nabla_x$.  Thus,  the gradient with respect to reference coordinates $\zeta$ is given by $\nabla_{\zeta} \equiv (\frac{\partial}{\partial \zeta_1}, \frac{\partial}{\partial \zeta_2}, \frac{\partial}{\partial \zeta_3}) = {\mathcal H} {\mathcal A} \nabla_x$.   So,  system (\ref{vv}) can be finally rewritten for $(\zeta , t)$ in $\hat\Omega\times(0,\infty)$ as:
\begin{equation} \label{vv2}
\left\{ \begin{array}{l}
\displaystyle{\frac{\partial \vec{u} }{\partial t}-\mu \nabla \cdot ({\mathcal A}^T {\mathcal H}^2 {\mathcal A} \nabla \vec{u})  + {\mathcal H} {\mathcal A} \nabla q} =  \vec{f} \circ \eta \circ \tau,
\vspace*{0.15cm} \\
\displaystyle{{\mathcal H} {\mathcal A} \nabla \cdot \vec{u} = 0} ,
\end{array}\right.
\end{equation}
where the second equation is equivalent to $ \nabla \cdot ({\mathcal A}^T {\mathcal H} \vec{u}) = 0$.

Throughout this section, as an abuse of notation,  expressions of the form $\vec{f} \circ \eta \circ \tau$ denote the compositions given by:
$$(\vec{f} \circ \eta \circ \tau) (\zeta,t) \equiv \vec{f}(\eta(\tau(\zeta),t),t) $$

Moreover, we have the estimate \cite{exi4}:
\begin{equation}  \nonumber
\| \vec{u} \|_{[W_p^{2,1}(\hat\Omega\times(0,T))]^3} + \| q \|_{ W_p^{1,0}(\hat\Omega\times(0,T))} \leq C(T) \sup \| \vec{F} \| 
\end{equation}
that,  from the results in Remark \ref{R3},  turns into:
\begin{equation}  \label{estv}
\| \vec{u} \|_{[W_p^{2,1}(\hat\Omega\times(0,T))]^3} + \| q \|_{ W_p^{1,0}(\hat\Omega\times(0,T))} \leq C(T) F \omega_{max}^2 \rho^2 . 
\end{equation}

\subsection{Estimates for the biological model}

Arguing in a similar way for the biological system (\ref{A}),  we obtain,  for instance,  for its first equation:
\begin{equation} \label{AA2} \hspace*{-.7cm}
\frac{\partial B }{\partial t} -\mu_A \nabla \cdot ({\mathcal A}^T {\mathcal H}^2 {\mathcal A} \nabla B)  = \big( (L \, \frac{P_1}{K_P+P_1} \,\frac{N_1+N_3}{K_N+N_1+N_3}) \circ \eta \circ \tau - (\gamma + \beta)  \big) B, 
\end{equation}
where $B$ represents algal concentration $A$ in reference coordinates $\zeta$, that is, $B=A\circ \eta \circ \tau$:
$$B(\zeta,t) = A(\eta(\tau(\zeta),t),t) = A(x,t).$$
(This process can be made in an analogous way for the rest of concentrations in system (\ref{A})).

Now, taking into account that $ \vec u \in [L^{p}((0,T);W^{2,p}(\hat\Omega))]^3$ with $p>3$,  we have from equation (\ref{eta}) that $\eta \in [W^{1,p}((0,T);W^{2,p}(\hat\Omega))]^3$.
Consequently,  from the definition of $\mathcal A$, we have that $\mathcal A  \in  [W^{1,p}(I;W^{1,p}(\hat\Omega))]^{3\times 3} $ and,  since $p>3$,  $\mathcal A  \in  [{\mathcal C}([0,T];{\mathcal C}(\overline{\hat\Omega})) ]^{3\times 3} $. Then, since from its definition the diagonal matrix ${\mathcal H} \in [{\mathcal C}^{\infty}(\overline{\hat\Omega})]^{3\times 3}$,  we obtain that the product ${\mathcal A}^T {\mathcal H}^2 {\mathcal A} \in  [W^{1,p}(I;W^{1,p}(\hat\Omega))]^{3\times 3} \subset [{\mathcal C}([0,T];{\mathcal C}(\overline{\hat\Omega})) ]^{3\times 3} $. 

On the other hand, all the second members in the biological system (\ref{A}) have first degree polynomial growth with respect to the different species concentations,  and all the initial conditions are bounded.  Thus,  all the assumptions of Theorem 1 of \cite{Bot} are fulfilled and, consequently, we obtain that the system  (\ref{A}) in reference coordinates has a unique nonnegative solution. In particular, for the algal concentration $B$ we have that:
\begin{equation}  \label{estA}
\begin{array}{l}
B \in L^{2}(I; H^1(\hat\Omega)) \cap L^{\infty}(I; L^{\infty}(\hat\Omega)) \cap \mathcal{C}([0,T]; L^2(\hat\Omega)),  \vspace*{0.1cm} \\
\| B \|_{L^{2}(I; H^1(\hat\Omega))} + \| B \|_{ L^{\infty}(I; L^{\infty}(\hat\Omega))} \leq \hat C(T) ,  \vspace*{0.15cm}\\
0 \leq B(\zeta,t) \leq \hat C(T),   \ \forall t \in I,  \ \xi \in  \hat\Omega ,
\end{array}
\end{equation}
(and similarly for the rest of concentrations in reference coordinates).

\begin{remark}\label{R3bis}
We must remark here that above wellposedness results are con\-gruent with the existence and regularity results previously obtained by the authors in \cite{Al2} for the biological system in Euler coordinates.
\end{remark}

\subsection{Existence of optimal solution}

Now, using above {\it a priori} estimates (\ref{estv}) and (\ref{estA}),  we will be able to demonstrate the existence of, at least, one solution of the optimal control problem $(\mathcal{P})$. 

From the boundedness of $B$ obtained in (\ref{estA}), and bearing in mind that $\det({\mathcal H} {\mathcal A}) = \det({\mathcal H}) \det( {\mathcal A}) = H$,  we deduce that cost functional
$$J(H,\omega) = - \int_{\Omega(T)} A(x,T) \, dx = - \int_{\hat\Omega} H \, B(\zeta,T) \, d\zeta $$
is bounded from below.  Then, there exists an infimum 
$$ d \equiv \inf_{\scriptsize \begin{array}{c} H_{min} \leq H \leq H_{max} \\ \omega_{min} \leq \omega \leq \omega_{max} \end{array}} J(H, \omega) $$
satisfying the state constraints,  and a minimizing sequence $\{ (H_n, \omega_n)\}_{n=1}^{\infty}\subset [H_{min}, H_{max}]\times [\omega_{min}, \omega_{max}]$ such that $J(H_n, \omega_n) \to d$,  as $n \to \infty$.

Moreover,  from the compactness of $[H_{min}, H_{max}] \times [\omega_{min}, \omega_{max}]\subset \mathbb{R}^2$, we obtain the existence of a subsequence (still denoted in the same way) such that $(H_n, \omega_n ) \to (\bar H, \bar\omega) \in [H_{min}, H_{max}]\times [\omega_{min}, \omega_{max}]$. 
Then,  from the definition of $\vec{f} = \vec{F}$ given in Subsection \ref{sub1}, we deduce that:
\begin{equation} \label{conf}
\vec{f}_n \equiv \vec{F}(H_n, \omega_n) \to \vec{\bar f} \equiv \vec{F}(\bar H, \bar\omega) \quad \mbox{in } L^{\infty}(I;L^{\infty}(\hat\Omega)).
\end{equation}

Let $(\vec{u}_n, q_n, B_n)$ solutions of the state systems associated to $(H_n, \omega_n)$. Then, from the boundedness obtained from estimates (\ref{estv}) and (\ref{estA}),  we have (maybe up to a subsequence) the following weak convergences to $(\vec{\bar u}, \bar q, \bar B)$:
\begin{equation}  \label{conA}
\begin{array}{l}
\vec{u}_n  \rightharpoonup \vec{\bar u}  \quad \mbox{in } [ L^{p}(I;W^{2,p}(\hat\Omega)) \cap W^{1,p}(I;L^{p}(\hat\Omega))]^3 \\
q_n  \rightharpoonup \bar q  \quad \mbox{in } L^{p}(I;W^{1,p}(\hat\Omega)) \\
B_n  \rightharpoonup \bar B  \quad \mbox{in } L^{2}(I;H^{1}(\hat\Omega)) \\
B_n  \rightharpoonup^{\ast} \bar B  \quad \mbox{in } L^{\infty}(I;L^{\infty}(\hat\Omega)) \\
B_n(T)  \rightharpoonup \bar B(T)  \quad \mbox{in } L^{2}(\hat\Omega) 
\end{array}
\end{equation}

Now,  from the characterization (\ref{eta}) of $\eta$ and the definition of matrix ${\mathcal A}$ we have that:
\begin{equation}  \label{coneta}
\begin{array}{l}
\eta_n  \rightharpoonup \bar\eta  \quad \mbox{in } [W^{1,p}(I;W^{2,p}(\hat\Omega)) \cap W^{2,p}(I;L^{p}(\hat\Omega))]^3 \\
{\mathcal A}_n  \rightharpoonup \bar{\mathcal A}  \quad \mbox{in } [W^{1,p}(I;W^{1,p}(\hat\Omega))]^{3\times 3} 
\end{array}
\end{equation}
that, from the Sobolev embedding theorem (since $p>3$),  implies the strong convergences:
\begin{equation}  \label{coneta2}
\begin{array}{l}
\eta_n  \to \bar\eta  \quad \mbox{in } [{\mathcal C}([0,T];{\mathcal C}^1(\overline{\hat\Omega}))]^3 \\
{\mathcal A}_n  \to \bar{\mathcal A}  \quad \mbox{in } [{\mathcal C}([0,T];{\mathcal C}(\overline{\hat\Omega})) ]^{3\times 3}
\end{array}
\end{equation}
and subsequently in $L^{\infty}(I;L^{\infty}(\hat\Omega))$. 

Similarly, from the definition of $\tau$,  we also have that:
\begin{equation}  \label{coneta3}
\begin{array}{l}
\tau_n  \to \bar\tau  \quad \mbox{in } [{\mathcal C}^{\infty}(\overline{\hat\Omega})]^3 \\
{\mathcal H}_n  \to \bar{\mathcal H}  \quad \mbox{in } [{\mathcal C}^{\infty}(\overline{\hat\Omega})]^{3\times 3} .
\end{array}
\end{equation}

Thus,  as a consequence of the last convergence in (\ref{conA}) and above conver\-gences,  we deduce that:
\begin{equation}  \nonumber %\label{conJ}
J(H_n,\omega_n) \equiv - \int_{\hat\Omega} H_n  B_n(\zeta,T) \, d\zeta \to - \int_{\hat\Omega} \bar H \bar B(\zeta,T) \, d\zeta \equiv J(\bar H, \bar\omega)=d.
\end{equation}

At this point, we need to assume more additional regularity,  mainly that 
$$\mbox{ (H1) }  \qquad \qquad {\mathcal A}_n  \to \bar{\mathcal A} \quad \mbox{in } [L^{\infty}(I;W^{1,\infty}(\hat\Omega))]^{3\times 3}. $$
Under this new assumption we can now deduce that, from the strong conver\-gences in $L^{\infty}$  and the weak convergences in $L^{p}$ given in (\ref{conf})-(\ref{coneta3}),  the equation
$$\frac{\partial \vec{u}_n }{\partial t}-\mu \nabla \cdot ({\mathcal A}_n^T {\mathcal H}_n^2 {\mathcal A}_n \nabla \vec{u}_n)  + {\mathcal H}_n {\mathcal A}_n \nabla q_n =  \vec{f}_n \circ \eta_n \circ \tau_n $$
converges weakly in $[L^{p}(I;L^{p}(\hat\Omega))]^3$ to
$$\frac{\partial \vec{\bar u} }{\partial t}-\mu \nabla \cdot (\bar{\mathcal A}^T \bar{\mathcal H}^2 \bar{\mathcal A} \nabla \vec{\bar u})  + \bar{\mathcal H} \bar{\mathcal A} \nabla \bar q =  \vec{\bar f} \circ \bar\eta \circ \bar\tau ,$$
and that the equation
$${\mathcal H}_n {\mathcal A}_n \nabla \cdot \vec{u}_n = 0 $$
converges weakly in $L^{p}(I;W^{1,p}(\hat\Omega))$ to
$$\bar{\mathcal H} \bar{\mathcal A} \nabla \cdot \vec{\bar u} = 0 ,$$
that is,  $(\vec{\bar u},  \bar q)$ is a solution of system (\ref{vv2}).

Similarly,  assuming more additional hypotheses on the regularity of the concentrations of the biological system (\ref{A}), mainly for $B=A\circ \eta \circ \tau$:
$$\mbox{ (H2) }  \qquad \qquad B  \mbox{ is bounded in }  L^{2}(I;H^{2}(\hat\Omega)) \cap H^{1}(I;L^{2}(\hat\Omega)), $$
we can obtain the weak convergence in $L^{2}(I;L^{2}(\hat\Omega))$ of equation (\ref{AA2}) -and the rest of equations of the biological system (\ref{A})- assuring that $\bar B$ corresponds to a solution of system (\ref{A}).

Finally,  arguing as above,  we can also obtain that state constraints (\ref{c1}) and (\ref{c2}) are satisfied.

Thus,  we obtain that:
$$ J(\bar H, \bar\omega) = \lim_{n\to\infty} J(H_n, \omega_n) = d = \inf_{\scriptsize \begin{array}{c} H_{min} \leq H \leq H_{max} \\ \omega_{min} \leq \omega \leq \omega_{max} \end{array}} J(H, \omega) \leq J(\bar H, \bar\omega) ,$$
which means that 
$$\min_{\scriptsize \begin{array}{c} H_{min} \leq H \leq H_{max} \\ \omega_{min} \leq \omega \leq \omega_{max} \end{array}} J(H, \omega) = J(\bar H, \bar\omega).$$

Summarizing,  we have demonstrated the following existence result:

\begin{theorem}
Assuming that above additional hypotheses (H1) and (H2) are satisfied,  then $(\bar H, \bar\omega)$ is a solution of the optimal control problem $(\mathcal{P})$.
\end{theorem}

\begin{remark}\label{R4}
It is worthwhile remarking here that we have only proved the existence of, at least, a solution of the optimal control problem.  However,  uniqueness of solution is not expected,  due to the high nonlinearity of the problem. 
On the other hand,  with respect to hypotheses (H1) and (H2),  these assumptions could be avoided, for instance,  by imposing more regularity for the second member $\vec{F}$ of the hydrodynamic system (\ref{v}), but this would not be a realistic supposition for the real-world problem we are dealing with.
\end{remark}

\section {Numerical implementation}

For the numerical resolution of above optimal control problem, we need first to give a suitable variational formulation of our state systems (\ref{v}) and (\ref{A}), within the Arbitrary Lagrangian Eulerian (ALE) framework. 
In this case, we set a standard variational formulation of both systems, employing the $\sigma$-transformation in order to deal with the free surface. (As it was demonstrated in \cite{Dec}, the classical $\sigma$-transformation is equivalent to the ALE formulation for a particular type of ALE mapping).

Since we are interested in employing the open-source module TELEMAC-3D \cite{Her} in our numerical computations,  we will use a space-time discretization,  with a time semi-discretization given by a suitable time step $\Delta t = T / N$ (which defines a set of discrete times $t_n=n \Delta t,  \,  n=0,1,\dots,N$),  and where the chosen spatial semi-discretization is based on a finite element method using prisms (whose 3D meshes can be easily varied with the $\sigma$-transformation as the free surface evolves). Moreover,  as an alternative to the classical method of characteristics, the advective terms appearing in the state system will be treated with the Multidimensional Upwind Residual Distribution (MURD) method.

Once we have computed a discrete approximation of the state variables (in particular,  the discretized velocity $\vec{v}_h^n(\cdot) \simeq \vec{v}(\cdot,t_n)$,  the discretized concentra\-tion of algae $A_h^n(\cdot) \simeq A(\cdot,t_n)$,  and the discretized concentration of dissolved oxygen $O_h^n(\cdot) \simeq O(\cdot,t_n)$,  for  $n=0,1,\dots,N$),  in order to calculate a discrete approximation $J_h$ of the cost functional $J$ given by (\ref{J}),  we can use any standard quadrature rule for numerical integration over the spatial domains:
\begin{equation}  \nonumber
J_h(H,\omega) = - \int_{\Omega_h^N} A_h^N(x) \, dx .
\end{equation}  
Finally,  we must also bear in mind that the imposed state constraints $(\ref{c1})-(\ref{c2})$ need to be incorporated into the approximated cost functional as added penalty terms.
So, we arrive to the discrete penalized cost function:
\begin{eqnarray*} 
&& \tilde J_h(H,\omega) = - \int_{\Omega_h^N} A_h^N(x) \, dx  + M_1 \, \max \Big\{ C_1 - \sum_{n=1}^N  \Delta t  \int_{\Omega_h^n} \| \vec{v}_h^n(x) \| \, dx , 0 \Big\} \nonumber \\
&& \qquad \qquad + M_2 \, \max \Big\{ C_2 - \min_{1 \leq n \leq N} \big\{ \int_{\Omega_h^n} O_h^n(x) \, dx \big\} , 0 \Big\}  
\label{JJ}
\end{eqnarray*}
where $M_1,  M_2 >0$ are penalty weights associated to state constraints (\ref {c1}) and (\ref {c2}),  respectively.

In this way, we will arrive to the discrete, bound-constrained minimization problem:
\begin{equation}  \nonumber
\min_{\scriptsize \begin{array}{c} H_{min} \leq H \leq H_{max} \\ \omega_{min} \leq \omega \leq \omega_{max} \end{array}} \tilde J_h(H,\omega) ,
\end{equation}  
whose solution can be obtained by any numerical optimization algorithm, in particular, and for the sake of simplicity,  by any derivative-free algorithm.  In the present case we will propose, for instance, the Nelder-Mead method \cite{7},  after the inclusion of a suitable penalty term to deal with the bound constraints on the control variables $(H,\omega)$, since Nelder-Mead algorithm is a method for unconstrained minimization.

\section {A computational example} \label{S4}

For the sake of conciseness,  we will present here only one numerical test from the many developed by the authors.   Our real-world scenario is a raceway with following dimensions (in meters): length of straight channels $L=20.0$,  width $W=2.0$,  and radii $r=0.2$ and $R=2.2$.  The time interval corresponds to one day, that is,  $T=86400$ seconds.  Numerical example presented here was obtained with the open-source hydrodynamics module TELEMAC-3D (compared and validated with those achieved with commercial program MIKE21). 

For the rotating paddlewheel,  we consider a force magnitude of $F=10.0$,  and paddles of length $\rho=0.4$,  and we fix, for technical reasons,  the coordinate $x_3^0=\rho+0.1=0.5$ (so that the paddle does not pass too close to the bottom of the raceway).  Coordinate $x_1^0$ is taken as $5.0$.  Moreover, due to pond symmetry,  we can restrict our study to only one of the two parallel straight channels -say the left one- and we fix the coordinate $x_2^0$ as the value corresponding to the central width of this left half of the raceway: $x_2^0=1.2$.  So,  the region of influence of the paddles for the force term $\vec{F}$ in (\ref{v}) is given by the horizontal cylindrical segment:
$$R(t) = \{ (x_1, x_2,x_3) \in \Omega(t) \ /  \  r \leq x_2 \leq R, \  (x_1-x_1^0)^2 + (x_3-x_3^0)^2 \leq \rho^2 \}, $$ which is independent of angular speed $\omega$,  depending only on height $H$.

For the other data of the state systems,  with time $t$ measured in seconds,  our future intention is to take temperature $\theta$ (in $^{\circ}$C) given by
$\theta(t) = \theta_0 + 2 \sin(2 \pi t / 86400)$, 
(this is,  considering a reference temperature $\theta_0 = 20$,  $\theta$ oscillates between 18 and 22 along the whole day),  and
incident light intensity $i$ given by expression $i(t) = \max \{ 0,  \sin(2 \pi t / 86400) \}$, 
which means that $i$ is null overnight. 
However,  here,  for the sake of simplicity,  we will take them as constant:  $\theta = 20$,  $i = 1$.
Finally,  for algal coefficients,  we use the physical parameters of {\it Chlorella} species \cite{Vis}. % (or {\it Scenedesmus} sp.)

With respect to the objective function,  we consider the state constraints thresholds $C_1=0.0$ and $C_2=4.0$ (that is,  constraint (\ref{c1}) on velocity is not taken into account),  and the control bounds $H_{min}=0.2$,  $H_{max}=0.5$,  $\omega_{min}=0.1$,  and $\omega_{max}=0.9$.

\begin{figure}[t]
\centering
\includegraphics[scale=0.45]{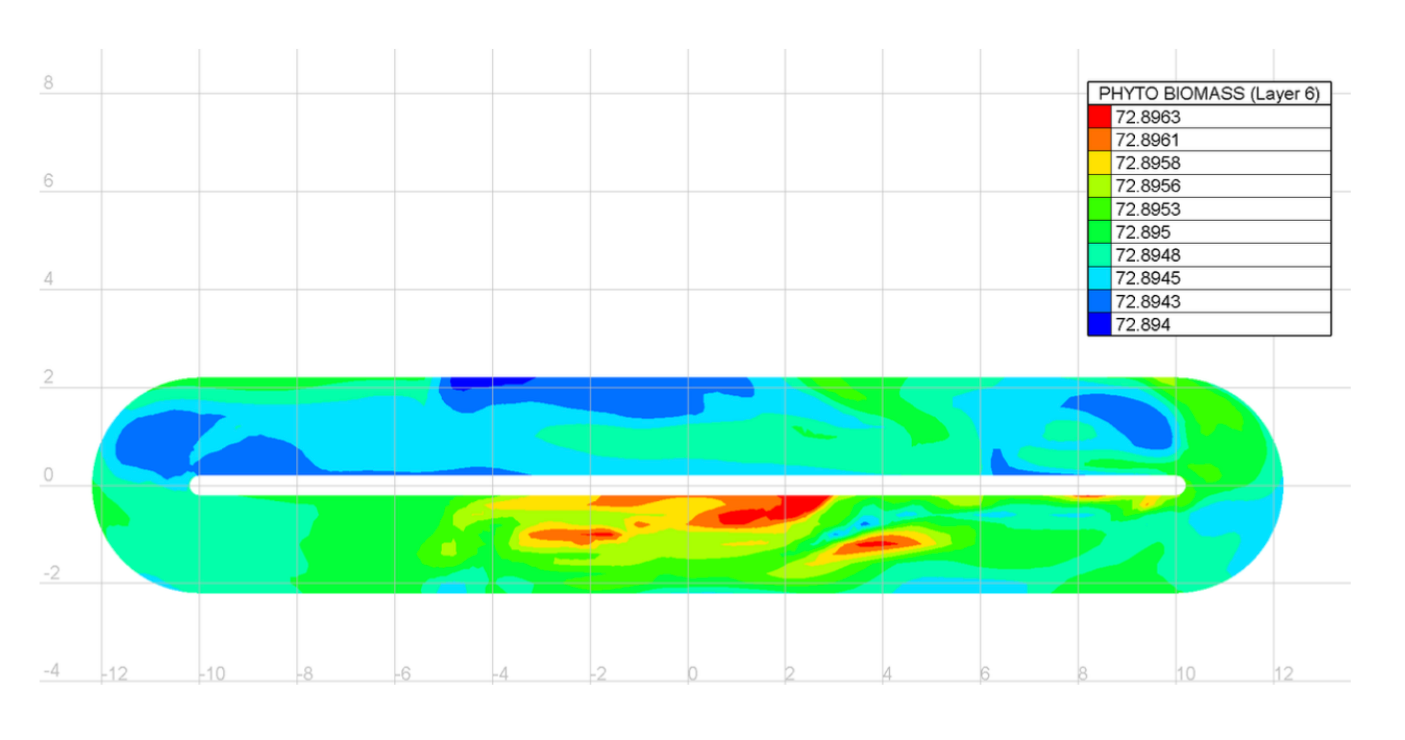}
\caption{Picture for algae concentration $A^h_N$ at final time $T=86400$ in the raceway central layer (corresponding to an initial algae concentration of $A_0=70.0$).} \label{figA}
\end{figure}

In the test shown here,  starting from random initial values $H=0.3$, $\omega=0.4$,  with a cost value $J_h(H,\omega)=-72.151$,  the Nelder-Mead algorithm arrives -after 41 iterations- to the optimal values $\bar H=0.2001$, $\bar\omega=0.4113$, corresponding to minimal cost value $J_h(\bar H,\bar\omega)=-72.895$.

\begin{figure}[t]
\centering
\includegraphics[width=11cm]{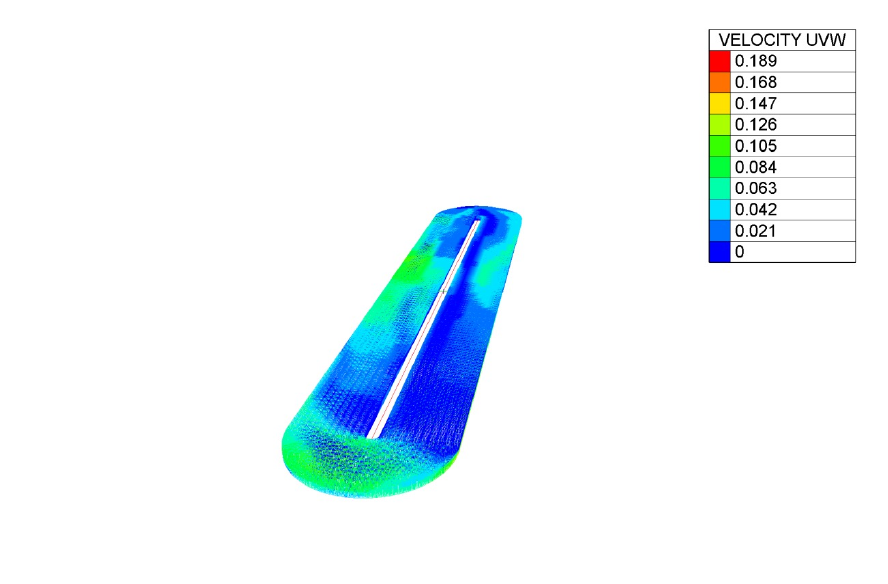}
\caption{Computational water velocities for the raceway mesh at final time $T$.} \label{figure5}
\end{figure}

In Fig.~\ref{figA} we show an example of numerical results corresponding to optimal algae concentration in the middle layer of the raceway at final time $T=86400$.   Fig.~\ref{figure5} shows the computational water velocities for the raceway mesh at final time $T$ for the optimized case.

To speed up the convergence of the optimization process,  alternative gradient-type methods will be proposed and tested in a forthcoming paper, where -if necessary- derivatives will be approximated by finite differences.

\section{Conclusions}

This paper introduces a novel strategy to optimize the performance of an open-channel raceway pond,  based on the application of the theory of optimal control for partial differential equations.  After analyzing the mathematical formulation of the optimal control problem,  we propose a full numerical algorithm for its computational resolution,  and present also some numerical tests for a real-world example.

Promising results achieved for this realistic case,  make manageable the inclusion of other design items to be optimized (for instance,  the optimal dimensions of the fishway,   the optimal location for the paddlewheel,  or the optimal length of the paddles),  and also the exploration of other alternative objectives and/or constraints in the cost functional to be minimized,  corres\-ponding to new interests of the stakeholders or to additional technological restrictions.

\section*{Acknowledgements}
This research was funded by Ministerio de Ciencia e Innovaci\'on (Spain) and NextGenerationEU (European Union) grant number TED2021-129324B-I00.  
The support provided by DHI with modelling system MIKE21 is also greatly appreciated.
The authors also thank the interesting suggestions of Carmen Rodr\'iguez, from the Department of Applied Mathematics in Universidade de
Santiago de Compostela (Spain).

\end{document}